\documentclass[letterpaper,conference,10pt]{IWMRRF}
\usepackage{amsmath,amssymb,amsfonts,amscd,amsbsy,amsxtra}
\usepackage{commath} 
\usepackage{graphicx}

\newcommand{\R}{\mathbb{R}}
\newcommand{\N}{\mathbb{N}}
\renewcommand{\epsilon}{\varepsilon}

\usepackage{url}       

\usepackage{amsmath}   
\usepackage{array}

\begin{document}

\title{An online slow manifold approach for efficient optimal control of multiple time-scale kinetics}

\author{\authorblockN{Marcus Heitel\authorrefmark{1},
Dirk Lebiedz\authorrefmark{1}}
\authorblockA{\authorrefmark{1}Institute for Numerical Mathematics, Ulm University, Germany}}

\maketitle

\begin{abstract}
  Chemical reactions modeled by ordinary differential equations are finite-dimensional dissipative dynamical systems
  with multiple time-scales. They are numerically hard to tackle -- especially when they enter an optimal
  control problem as ``infinite-dimensional'' constraints. Since discretization of such problems usually results in high-dimensional nonlinear problems, model (order)
  reduction via slow manifold computation seems to be an attractive approach. We discuss the use of slow
  manifold computation methods in order to solve optimal control problems more efficiently having real-time applications in view.
\end{abstract}

\section{Introduction}

Chemical kinetics with multiple time scales and their control involve highly stiff and often high-dimensional ordinary differential equations (ODE). 
This poses hard challenges to the numerical solution and is the reason why model reduction methods are considered. The dynamics can be simplified by focusing on the long time behavior of such systems (leaving fast transients unresolved) and
calculating fast modes as functions of the slow ones. Ideally this leads to low dimensional manifolds  in high-dimensional state space.
In the special case of singularly perturbed system, they are understood quite well and called slow invariant manifolds.

An open problem is how the slow manifolds can be used to simplify the solution of optimal control problems (OCP) that involve multiple time scale ODE constraints.

\section{Slow Manifold Computation}

In dissipative dynamical systems geometrically the bund\-ling of trajectories (on a fast time scale) to low-dimensional manifolds is observed. 
Once trajectories reach the neighborhood of the slow manifold, they will evolve slowly and will never leave this manifold neighborhood.
Thus, this manifold is called {\em slow invariant attracting manifold} (SIAM).

The aim of slow manifold computation techniques is to approximately compute the SIAM as the graph of a function of only a few selected species (so called 
{\em reaction progress variables}). Thus, manifold-based model reduction generate a function $h:\R^{n_s} \rightarrow \R^{n_f}$ 
($n_s$ is the number of slow variables resp. reaction progress variables and $n_f$ is the number of fast variables), 
such that $\left(z_s,h(z_s)\right)$ approximates points of the SIAM.

\vspace*{1em}
In order to investigate optimal control benchmark problems, we consider singularly perturbed systems, i.e. systems where the ODE can be transformed into the following 
form:
\begin{subequations}\label{formula:sps}
	\begin{alignat}{2}
	\dot{z_s}(t) &= f_s(z_s(t),z_f(t))\\
	\varepsilon \dot{z_f}(t) &= f_f(z_s(t),z_f(t)).
	\end{alignat}
\end{subequations}
Two methods relevant in our context for the approximative calculation of the SIAM are briefly reviewed in the following subsections.

\subsection{Zero Derivative Principle}

The main idea of the Zero Derivative Principle (ZDP) \cite{Gear2005},\cite{Zagaris2009} for model reduction of singularly perturbed systems is to identify for given values of the slow variables $z_s^*$ a point $z_f^*$ such that the higher-order time derivatives of fast components vanish, i.e 
\begin{align}
\frac{\dif{}^m f_f(z_s^*,z_f^*)}{\dif t^m} = 0 \qquad \text{for a given } m \in \N.
\end{align}

\subsection{Method of Lebiedz and Unger}

Another approach proposed by Lebiedz and Unger \cite{Lebiedz2016} is motivated geometrically: Among arbitrary trajectories of (\ref{formula:sps})
for which the slow components end within the time $t_1-t_0$ in the state $z_s^*$ the corresponding part of the trajectory on the SIAM is characterized by 
the smallest curvature (see also \cite{Lebiedz2004},\cite{Lebiedz2011a}). 
This motivates optimization problem (\ref{formula:Unger}) which is a variational boundary value problem (BVP).
\begin{subequations} \label{formula:Unger}
	\begin{alignat}{2}
	\min\limits_{z(\cdot)=\big(z_s(\cdot),z_f(\cdot)\big)} \quad & & & \hspace*{-5mm} \norm{\ddot{z}(t_0)}_2^2 \\
	\text{ s.t. } \hspace*{1cm}	&  &\dot z_s &= f_s\big(z_s,z_f\big), \quad t \in [t_0,t_1] \\
	& & \epsilon \dot z_f &= f_f\big(z_s,z_f\big), \quad t \in [t_0,t_1]\\
	&  &	z_s(t_1) &= z_s^*.
	\end{alignat}
\end{subequations}
In our application context we also use the local reformulation of problem (\ref{formula:Unger}), where $t_0 = t_1$.

\section{Optimal Control}
One of our research interests is to solve optimal control problems involving multiple time scales as it appears frequently
e.g. in the field of chemical engineering. Thus, we consider the following (typically high-dimensional) OCP:
\begin{subequations} \label{formula:origOCP}
	\begin{alignat}{2}
	& \quad \, \min_{z_s,z_f,u} & 		\int_{0}^{T}	&L(z_s,z_f,u) \dif t  \\
	&\mbox{ subject to } \quad	& 		   \dot{z_s}	&= f_s(z_s,z_f,u)\\
	&							& \epsilon \dot{z_f}	&= f_f(z_s,z_f,u)\\
	&							& z_s(0)  				&=z_s^{(0)}, \; z_f(0) = z_f^{(0)}
	\end{alignat}
\end{subequations} 

Applying the model reduction methods presented in the last section and assuming the control $u$ to be a slow
variable, yields the lower dimensional problem (cf. \cite{Rehberg2013})
\begin{subequations}\label{formula:redOCP}
	\begin{alignat}{2} 
	& \quad \, \min_{z_s,u} 	& 		\int_{0}^{T}	&L(z_s,h(z_s,u),u) \dif t  \\
	&\mbox{ subject to } \quad	& 		   \dot{z_s}	&= f_s(z_s,h(z_s,u),u) \label{formula:redOCP_ode} \\
	&							& z_s(0)  				&=z_s^{(0)}.
	\end{alignat}
\end{subequations} 
This systems has the advantage, that it has significantly less optimization variables and the ODE (\ref{formula:redOCP}) 
is less stiff, which makes it solvable by fast explicit numerical integrators compared to implicit methods required for stiff ODE.
However, numerical solution methods for OCPs like the multiple shooting method need repeated evaluation of the function $h$ as well as its 
partial derivatives $h_{z_s}$ and $h_u$. 

Therefore, it would be beneficial to combine the calculation of the SIM and the optimal control problem. This is obviously possible, if the approximation
$h(z_s)$ of the SIM can be formulated as a (nonlinear) root finding problem $r(z_s,z_f,u) = 0$, e.g. with the ZDP method. Thus, we propose to solve
the following OCP instead of (\ref{formula:redOCP}):
\begin{subequations} \label{formula:liftedOCP}
\begin{alignat}{2} 
& \quad \, \min_{z_s,z_f,u} & 		\int_{0}^{T}	&L(z_s,z_f,u) \dif t  \\
&\mbox{ subject to } \quad	& 		   \dot{z_s}	&= f_s(z_s,z_f,u)\\
&							& 0						&= r(z_s,z_f,u)\\
&							& z_s(0)  				&=z_s^{(0)}.
\end{alignat}
\end{subequations} 
\section{Application to Chemical Reactions}

We apply the ideas presented in the last sections to a benchmark OCP motivated by the Michaelis-Menten-Henri mechanism 
\begin{align}\label{formula:enzyme}
S + E \leftrightharpoons SE \rightarrow P + E,
\end{align}
modeling the reaction of substrate $S$ to a product $P$ via a substrate-enzyme-complex $SE$ with the help of enzyme $E$.
Simplifying the ODE  given by (\ref{formula:enzyme}) and introducing an artificial objective function yields 
OCP (\ref{formula:enzymeOCP}).
\begin{subequations}\label{formula:enzymeOCP}
	\begin{alignat}{1}
	\min_{z_s,z_f,u}\quad 	& \int_{0}^{5} -50z_f + u^2 \dif t \\
	\mbox{s.t.} \hspace*{2mm} \quad				&  
	\begin{array}{rrrll}
	\dot{z_s} 				&=&- z_s&+ \left(z_s + 0.5\right)z_f + u, 	\\
	\epsilon \dot{z_f} 		&=& z_s &- \left(z_s+1\right)z_f , \\
	z_s(0)					&=& 1, 	& &
	\end{array} 
	\end{alignat}
\end{subequations} 
where the control $u(t) \in [0,10]$ represents the possibility to add
some substrate (corresponds to variable $z_s$) to the system and $\epsilon$ describes the
time-scale separation (between the time evolution of $z_s$ and $z_f$).

\begin{figure}[t]
	\centering
	\includegraphics[scale=0.4]{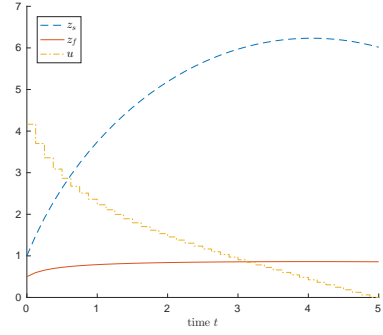}
	\caption{Numerical Solution of (\ref{formula:enzymeOCP}).} \label{fig:enzyme}
\end{figure}

Figure \ref{fig:enzyme} shows the results of the numerical solution of (\ref{formula:enzymeOCP}) using the multiple-shooting scheme
with an implicit Radau-2A integrator. If we refer to the solution of the proposed OCP (\ref{formula:liftedOCP}) as 
$(z_s^{\text{app}},z_f^{\text{app}},u^{\text{app}})$ and to the solution of (\ref{formula:enzymeOCP}) as 
$(z_s^{\text{orig}},z_f^{\text{orig}},u^{\text{orig}})$, then it holds
\begin{align}
	& \max\left\lbrace \lvert \lvert 	z_s^{\text{orig}}-z_s^{\text{app}} \rvert \rvert_{\infty}, 
	\lvert \lvert	z_f^{\text{orig}}-z_f^{\text{app}} \rvert \rvert_{\infty},
	\lvert \lvert	u^{\text{orig}}-u^{\text{app}} \rvert \rvert_{\infty} \right\rbrace \nonumber \\
	&= \lvert \lvert 	z_s^{\text{orig}}-z_s^{\text{app}} \rvert \rvert_{\infty} \approx 0.05,
\end{align} 
which gives a relative error of $\approx 0.2 \%$ for both objective functional value and $\lvert \lvert z_s^{\text{orig}}-z_s^{\text{app}} \rvert \rvert_{\infty} /\lvert \lvert z_s^{\text{orig}} \rvert \rvert_{\infty}$.
Although, the proposed method uses exactly as many variables than the original OCP, we observe a speed up of factor 4
for solving OCP (\ref{formula:enzymeOCP}) due to the use of an explicit integration scheme.





\end{document}